\documentclass[a4paper,12pt,oneside]{article}
\usepackage{amsmath,amsfonts,amssymb,amsmath,latexsym,amsthm,enumerate}

\newtheorem{thm}{Theorem}

\theoremstyle{definition}

\theoremstyle{plain}
\newtheorem{lem}[thm]{Lemma}

\newtheorem{prop}[thm]{Proposition}

\begin{document}
\title {Higher-order Bernoulli and poly-Bernoulli mixed type polynomials}
\author{by \\Dae San Kim and Taekyun Kim}\date{}\maketitle

\begin{abstract}
\noindent In this paper, we consider higher-order Bernoulli and poly-Bernoulli mixed type polynomials and we give some interesting identities of those polynomials arising from umbral calculus.
\end{abstract}

\section{Introduction}

The classical polylogarithmic function $Li_{k}(x)$ is
\begin{equation}\label{eq:1}
Li_{k}(x)=\sum_{n=1}^{\infty}\frac{x^{n}}{n^{k}},\,\,\,\,k\in\mathbf{Z},\,\,\,\,(\text{see}\,\, \lbrack 3,5\rbrack).
\end{equation}
The poly-Bernoulli polynomials are defined by the generating function to be
\begin{equation}\label{eq:2}
\frac{Li_{k}\left(1-e^{-t}\right)}{1-e^{-t}}e^{xt}=\sum_{n=0}^{\infty}B_{n}^{(k)}(x)\frac{t^{n}}{n!},\,\,\,(\text{see}\,\, \lbrack 3,5\rbrack),
\end{equation}
and the Bernoulli polynomials of order $r (r\in\mathbf{Z})$ are given by
\begin{equation}\label{eq:3}
\left(\frac{t}{e^{t}-1}\right)^{r}e^{xt}=\sum_{n=0}^{\infty}\mathbb{B}_{n}^{(r)}(x)\frac{t^{n}}{n!},\,\,\,(\text{see}\,\, \lbrack 2,4,7\rbrack).
\end{equation}
When $x=0$, $B_{n}^{(k)}=B_{n}^{(k)}(0)$ are called the poly-Bernoulli numbers and $\mathbb{B}_{n}^{(r)}=\mathbb{B}_{n}^{(r)}(0)$ are called the Bernoulli numbers of order $r$.
In the special case, $r=1$, $\mathbb{B}_{n}^{(1)}(x)=B_{n}(x)$ are called the Bernoulli polynomials. When $x=0$, $B_{n}=B_{n}(0)$ are called the ordinary Bernoulli numbers.\\
The higher-order Bernoulli and poly-Bernoulli mixed type polynomials are defined by the generating function to be
\begin{equation}\label{eq:4}
\left(\frac{t}{e^{t}-1}\right)^{r}\frac{Li_{k}\left(1-e^{-t}\right)}{1-e^{-t}}e^{xt}=\sum_{n=0}^{\infty}s_{n}^{(r,k)}(x)\frac{t^{n}}{n!},\,\,\,\,(\text{see}\,\, \lbrack 5\rbrack).
\end{equation}
From (\ref{eq:2}), (\ref{eq:3}) and (\ref{eq:4}), we note that
\begin{align}\label{eq:5}
s_{n}^{(r,k)}(x)&=\sum_{l=0}^{n}\binom{n}{l}B_{n-l}^{(k)}\mathbb{B}_{l}^{(r)}(x)\\
&=\sum_{l=0}^{n}\binom{n}{l}\mathbb{B}_{n-l}^{(r)}B_{l}^{(k)}(x).\nonumber
\end{align}
When $x=0$, $s_{n}^{(r,k)}=s_{n}^{(r,k)}(0)$ are called the higher-order Bernoulli and poly-Bernoulli mixed type numbers.\\
Let $\mathcal{F}$ be the set of all formal power series in variable $t$ over $\mathbf{C}$ with
\begin{equation}\label{eq:6}
\mathcal{F}=\left\{f(t)=\sum_{k=0}^{\infty}a_{k}\frac{t^{k}}{k!}\Bigg\vert a_{k}\in\mathbf{C}\right\}.
\end{equation}
Let $\mathbb{P}=\mathbf{C}\lbrack t\rbrack$ and let $\mathbb{P}^{*}$ be the vector space of all linear functional on $\mathbb{P}$. $\left\langle L\vert p(x)\right\rangle$ denotes the acition of linear functional $L$ on the polynomial $p(x)$, and it is well known that the vector space operations on $\mathbb{P}^{*}$ are defined by $\left\langle L+M\vert p(x)\right\rangle=\left\langle L\vert p(x)\right\rangle+\left\langle M\vert p(x)\right\rangle$, $\left\langle cL\vert p(x)\right\rangle=c\left\langle L\vert p(x)\right\rangle$, where $c$ is a complex constant.
For $f(t)\in\mathcal{F}$ with $f(t)=\sum_{k=0}^{\infty}a_{k}\frac{t^{k}}{k!}$, let us define the linear functional on $\mathbb{P}$ by setting
\begin{equation}\label{eq:7}
\left\langle f(t)\vert x^{n}\right\rangle=a_{n},\,\,\,\,(n\geq 0),\,\,\,\,(\text{see}\,\, \lbrack 8,9\rbrack).
\end{equation}
From (\ref{eq:6}) and (\ref{eq:7}), we note that
\begin{equation}\label{eq:8}
\left\langle t^{k}\big\vert x^{n}\right\rangle=n!\delta_{n,k},\,\,\,\,(n, k\geq 0),
\end{equation}
where $\delta_{n,k}$ is the Kronecker's symbol.\\
Let $f_{L}(t)=\sum_{k=0}^{\infty}\left\langle L\big\vert x^{k}\right\rangle\frac{t^{k}}{k!}$.
Then, by (\ref{eq:8}), we see that $\left\langle f_{L}(t)\big\vert x^{n}\right\rangle=\left\langle L\big\vert x^{n}\right\rangle$. Additionally, the map $L\longmapsto f_{L}(t)$ is a vector space isomorphism from $\mathbb{P}^{*}$ onto $\mathcal{F}$. Henceforth, $\mathcal{F}$ denotes both the algebra of the formal power series in $t$ and the vector space of all linear functionals on $\mathbb{P}$, and so an element $f(t)$ of $\mathcal{F}$ will be thought as both a formal power series and a linear functional. We call $\mathcal{F}$ the umbral algebra. The umbral calculus is the study of umbral algebra. The order $O\left(f(t) \right)$ of the power series $f(t)\neq0$ is the smallest integer for which $a_{k}$ does not vanish. If $O(f(t))=0$, then $f(t)$ is called an invertible series. If $O(f(t))=1$, then $f(t)$ is called a delta series. For $f(t)\in\mathcal{F}$ and $p(x)\in\mathbb{P}$, we have
\begin{equation}\label{eq:9}
f(t)=\sum_{k=0}^{\infty}\left\langle f(t)\big\vert x^{k}\right\rangle\frac{t^{k}}{k!},\quad p(x)=\sum_{k=0}^{\infty}\left\langle t^{k}\big\vert p(x)\right\rangle\frac{x^{k}}{k!}.
\end{equation}
Thus, by (\ref{eq:9}), we get
\begin{equation}\label{eq:10}
p^{(k)}(0)=\left\langle t^{k}\big\vert p(x)\right\rangle=\left\langle 1\big\vert p^{(k)}(x)\right\rangle,\,\,\,\,(\text{see}\,\, \lbrack 8,9\rbrack),
\end{equation}
where $p^{(k)}(x)=\frac{d^{k}p(x)}{dx^{k}}$.\\
From (\ref{eq:10}), we have
\begin{equation}\label{eq:11}
t^{k}p(x)=p^{(k)}(x)=\frac{d^{k}p(x)}{dx^{k}}.
\end{equation}
By (\ref{eq:11}), we easily see that
\begin{equation}\label{eq:12}
e^{yt}p(x)=p(x+y),\quad\left\langle e^{yt}\big\vert p(x)\right\rangle=p(y).
\end{equation}
For $f(t), g(t)\in\mathcal{F}$ with $O(f(t))=1$, $O(g(t))=0$, there exists a unique sequence $s_{n}(x)$ of polynomials such that $\left\langle g(t)f(t)^{k}\big\vert s_{n}(x)\right\rangle=n!\delta_{n,k}$, for $n, k\geq 0$. The sequence $s_{n}(x)$ is called the Sheffer sequence for $\left(g(t),f(t)\right)$, which is denoted by $s_{n}(x)\sim\left(g(t),f(t)\right)$.\\
Let $p(x)\in\mathbb{P}$ and $f(t)\in\mathcal{F}$. Then we see that
\begin{equation}\label{eq:13}
\left\langle f(t)\vert xp(x)\right\rangle=\left\langle \partial_{t}f(t)\vert p(x)\right\rangle=\left\langle f'(t)\vert p(x)\right\rangle,\,\,\,\,(\text{see}\,\,\lbrack 8\rbrack).
\end{equation}
For $s_{n}(x)\sim\left(g(t),f(t)\right)$, we have
\begin{equation}\label{eq:14}
\frac{1}{g\left(\bar{f}(t)\right)}e^{y\bar{f}(t)}=\sum_{k=0}^{\infty}s_{k}(y)\frac{t^{k}}{k!},\,\,\,\,\text{for all}\,\, y\in\mathbf{C},
\end{equation}
where $\bar{f}(t)$ is the compositional inverse for $f(t)$ with $\bar{f}\left(f(t)\right)=t$, and
\begin{equation}\label{eq:15}
f(t)s_{n}(x)=ns_{n-1}(x),\,\,\,\,(\text{see}\,\, \lbrack 8,9\rbrack),
\end{equation}
Let $s_{n}(x)\sim\left(g(t),t\right)$. Then we see that
\begin{equation}\label{eq:16}
s_{n+1}(x)=\left(x-\frac{g'(t)}{g(t)}\right)s_{n}(x),\,\,\,\,(\text{see}\,\,\lbrack 8\rbrack).
\end{equation}
For $s_{n}(x)\sim\left(g(t),f(t)\right)$, $r_{n}(x)\sim\left(h(t),l(t)\right)$, we have
\begin{equation}\label{eq:17}
s_{n}(x)=\sum_{m=0}^{n}c_{n,m}r_{m}(x),
\end{equation}
where
\begin{equation}\label{eq:18}
c_{n,m}=\frac{1}{m!}\left\langle \frac{h\left(\bar{f}(t)\right)}{g\left(\bar{f}(t)\right)}l\left(\bar{f}(t)\right)^{m}\Bigg\vert x^{n}\right\rangle,\,\,\,\,(\text{see}\,\,\lbrack 8,9\rbrack).
\end{equation}
In this paper, we study higher-order Bernoulli and poly-Bernoulli mixed type polynomials and we give some interesting identities of those polynomials arising from umbral calculus.

\section{Higher-order Bernoulli and poly-Bernoulli mixed type polynomials}

From (\ref{eq:4}) and (\ref{eq:14}), we note that
\begin{equation}\label{eq:19}
s_{n}^{(r,k)}(x)\sim\left(g_{r,k}(t)=\left(\frac{e^{t}-1}{t}\right)^{r}\frac{1-e^{-t}}{Li_{k}\left(1-e^{-t}\right)},t\right).
\end{equation}
Thus, by (\ref{eq:15}), we get
\begin{equation}\label{eq:20}
ts_{n}^{(r,k)}(x)=\frac{d}{dx}s_{n}^{(r,k)}(x)=ns_{n-1}^{(r,k)}(x).
\end{equation}
From (\ref{eq:4}), we have
\begin{equation}\label{eq:21}
s_{n}^{(r,k)}(x)=\sum_{l=0}^{n}\binom{n}{l}s_{l}^{(r,k)}x^{n-l}=\sum_{l=0}^{n}\binom{n}{l}s_{n-l}^{(r,k)}x^{l}.
\end{equation}
First, we observe that
\begin{equation}\label{eq:22}
s_{n}^{(r,k)}(x)=\frac{1}{g_{r,k}(t)}x^{n}=\left(\frac{t}{e^{t}-1}\right)^{r}\left(\frac{Li_{k}\left(1-e^{-t}\right)}{1-e^{-t}}\right)x^{n}.
\end{equation}
In $\lbrack 3\rbrack$, it is known that
\begin{equation}\label{eq:23}
\frac{Li_{k}\left(1-e^{-t}\right)}{1-e^{-t}}x^{n}=\sum_{m=0}^{n}\frac{1}{(m+1)^{k}}\sum_{j=0}^{m}(-1)^{j}\binom{m}{j}(x-j)^{n}.
\end{equation}
Thus, by (\ref{eq:22}) and (\ref{eq:23}), we get
\begin{align}\label{eq:24}
s_{n}^{(r,k)}(x)&=\sum_{m=0}^{n}\frac{1}{(m+1)^{k}}\sum_{j=0}^{m}(-1)^{j}\binom{m}{j}\left(\frac{t}{e^{t}-1}\right)^{r}(x-j)^{n}\\
&=\sum_{m=0}^{n}\frac{1}{(m+1)^{k}}\sum_{j=0}^{m}(-1)^{j}\binom{m}{j}\mathbb{B}_{n}^{(r)}(x-j).\nonumber
\end{align}
Therefore, by (\ref{eq:24}), we obtain the following proposition

\begin{prop}\label{eq:prop1}
For $n\in\mathbf{Z}_{\geq 0}$, $r, k\in\mathbf{Z}$, we have
\begin{equation*}
s_{n}^{(r,k)}(x)=\sum_{m=0}^{n}\frac{1}{(m+1)^{k}}\sum_{j=0}^{m}(-1)^{j}\binom{m}{j}\mathbb{B}_{n}^{(r)}(x-j).
\end{equation*}
\end{prop}

\noindent From (\ref{eq:3}), we can easily derive the following equation:
\begin{equation}\label{eq:25}
\mathbb{B}_{n}^{(r)}(x)=\sum_{l=0}^{n}\binom{n}{l}\mathbb{B}_{n-l}^{(r)}x^{l}.
\end{equation}
By (\ref{eq:24}) and (\ref{eq:25}), we get
\begin{equation}\label{eq:26}
s_{n}^{(r,k)}(x)=\sum_{l=0}^{n}\left\{\binom{n}{l}\mathbb{B}_{n-l}^{(r)}\sum_{m=0}^{n}\frac{1}{(m+1)^{k}}\sum_{j=0}^{m}(-1)^{j}\binom{m}{j}\right\}(x-j)^{l}.
\end{equation}
In $\lbrack 3\rbrack$, it is known that
\begin{equation}\label{eq:27}
\frac{Li_{k}\left(1-e^{-t}\right)}{1-e^{-t}}x^{n}=\sum_{j=0}^{n}\left\{\sum_{m=0}^{n-j}\frac{(-1)^{n-m-j}}{(m+1)^{k}}\binom{n}{j}m!S_{2}(n-j,m)\right\}x^{j},
\end{equation}
where $S_{2}(n,m)$ is the Stirling number of the second kind.\\
From (\ref{eq:22}) and (\ref{eq:27}), we have
\begin{align}\label{eq:28}
s_{n}^{(r,k)}(x)&=\sum_{j=0}^{n}\left\{\sum_{m=0}^{n-j}\frac{(-1)^{n-m-j}}{(m+1)^{k}}\binom{n}{j}m!S_{2}(n-j,m)\right\}\left(\frac{t}{e^{t}-1}\right)^{r}x^{j}\\
&=\sum_{j=0}^{n}\left\{\sum_{m=0}^{n-j}\frac{(-1)^{n-m-j}}{(m+1)^{k}}\binom{n}{j}m!S_{2}(n-j,m)\right\}\mathbb{B}_{j}^{(r)}(x)\nonumber\\
&=\sum_{l=0}^{n}\left\{\sum_{j=l}^{n}\sum_{m=0}^{n-j}(-1)^{n-m-j}\binom{n}{j}\binom{j}{l}\frac{m!}{(m+1)^{k}}S_{2}(n-j,m)\mathbb{B}_{j-l}^{(r)}\right\}x^{l}.\nonumber
\end{align}
From (\ref{eq:16}) and (\ref{eq:19}), we have
\begin{equation}\label{eq:29}
s_{n+1}^{(r,k)}(x)=\left(x-\frac{g_{r,k}'(t)}{g_{r,k}(t)}\right)s_{n}^{(r,k)}(x),
\end{equation}
where
\begin{align}\label{eq:30}
\frac{g_{r,k}'(t)}{g_{r,k}(t)}&=\left(\log{g_{r,k}(t)}\right)'\\
&=\left(r\log{\left(e^{t}-1\right)}-r\log{t}+\log{\left(1-e^{-t}\right)}-\log{Li_{k}\left(1-e^{-t}\right)}\right)'\nonumber\\
&=\frac{rte^{t}-re^{t}+r}{t\left(e^{t}-1\right)}+\frac{t}{e^{t}-1}\left(\frac{Li_{k}\left(1-e^{-t}\right)-Li_{k-1}\left(1-e^{-t}\right)}{tLi_{k}\left(1-e^{-t}\right)}\right).\nonumber
\end{align}
By (\ref{eq:29}) and (\ref{eq:30}), we get
\begin{align}\label{eq:31}
s_{n+1}^{(r,k)}(x)&=xs_{n}^{(r,k)}(x)-\left(\frac{t}{e^{t}-1}\right)^{r}\frac{Li_{k}\left(1-e^{-t}\right)}{1-e^{-t}}\left(\frac{rte^{t}-re^{t}+r}{t(e^{t}-1)}\right)x^{n}\\
&\quad -\left(\frac{t}{e^{t}-1}\right)^{r+1}\frac{Li_{k}\left(1-e^{-t}\right)-Li_{k-1}\left(1-e^{-t}\right)}{t\left(1-e^{-t}\right)}x^{n}.\nonumber
\end{align}
It is easy to show that
\begin{equation}\label{eq:32}
\frac{rte^{t}-re^{t}+r}{e^{t}-1}=\frac{r}{2}t+\cdots,\\
\frac{Li_{k}\left(1-e^{-t}\right)-Li_{k-1}\left(1-e^{-t}\right)}{1-e^{-t}}=\left(\frac{1}{2^{k}}-\frac{1}{2^{k-1}}\right)t+\cdots
\end{equation}
For any formal power series $f(t)$ with $O\left(f(t)\right)\geq 1$, we have
\begin{equation}\label{eq:33}
\frac{f(t)}{t}x^{n}=\frac{f(t)}{t}\frac{1}{n+1}tx^{n+1}=\frac{1}{n+1}f(t)x^{n+1}.
\end{equation}
By (\ref{eq:31}), (\ref{eq:32}) and (\ref{eq:33}), we get
\begin{align}\label{eq:34}
s_{n+1}^{(r,k)}(x)&=xs_{n}^{(r,k)}(x)-\frac{r}{n+1}\sum_{l=0}^{n}\binom{n+1}{l}(-1)^{n+1-l}B_{n+1-l}s_{l}^{(r,k)}(x)\\
&\quad -\frac{1}{n+1}\left\{s_{n+1}^{(r+1,k)}(x)-s_{n+1}^{(r+1,k-1)}(x)\right\}.\nonumber
\end{align}
Therefore, by (\ref{eq:34}), we obtain the following theorem.

\begin{thm}\label{eq:thm2}
For $r,k\in\mathbf{Z}$ and $n\in\mathbf{Z}_{\geq 0}$, we have
\begin{align*}
s_{n+1}^{(r,k)}(x)&=xs_{n}^{(r,k)}(x)-\frac{r}{n+1}\sum_{l=0}^{n}\binom{n+1}{l}(-1)^{n+1-l}B_{n+1-l}s_{l}^{(r,k)}(x)\\
&\quad -\frac{1}{n+1}\left\{s_{n+1}^{(r+1,k)}(x)-s_{n+1}^{(r+1,k-1)}(x)\right\}.
\end{align*}
\end{thm}

\noindent From (\ref{eq:5}), we have
\begin{align}\label{eq:35}
txs_{n}^{(r,k)}(x)&=\sum_{l=0}^{n}\binom{n}{l}\mathbb{B}_{n-l}^{(r)}t\left(xB_{l}^{(k)}(x)\right)\\
&=\sum_{l=0}^{n}\binom{n}{l}\mathbb{B}_{n-l}^{(r)}\left\{lxB_{l-1}^{(k)}(x)+B_{l}^{(k)}(x)\right\}\nonumber\\
&=nx\sum_{l=0}^{n-1}\binom{n-1}{l}\mathbb{B}_{n-1-l}^{(r)}B_{l}^{(k)}(x)+\sum_{l=0}^{n}\binom{n}{l}\mathbb{B}_{n-l}^{(r)}B_{l}^{(k)}(x)\nonumber\\
&=nxs_{n-1}^{(r,k)}(x)+s_{n}^{(r,k)}(x).\nonumber
\end{align}
It is easy to show that
\begin{equation}\label{eq:36}
s_{n}^{(r,k)}(x)=\left(\frac{t}{e^{t}-1}\right)^{r}\frac{Li_{k}\left(1-e^{-t}\right)}{1-e^{-t}}x^{n}=\left(\frac{t}{e^{t}-1}\right)^{r}B_{n}^{(k)}(x).
\end{equation}
Applying $t$ on the both sides of (\ref{eq:22}) and using (\ref{eq:36}), we get
\begin{align}\label{eq:37}
&(n+1)s_{n}^{(r,k)}(x)\\
&=nxs_{n-1}^{(r,k)}(x)+s_{n}^{(r,k)}(x)-\frac{r}{n+1}\sum_{l=1}^{n}\binom{n+1}{l}(-1)^{n+1-l}B_{n+1-l}ls_{l-1}^{(r,k)}(x)\nonumber\\
&\quad -\frac{1}{n+1}\left\{(n+1)s_{n}^{(r+1,k)}(x)-(n+1)s_{n}^{(r+1,k-1)}(x)\right\}\nonumber\\
&=nxs_{n-1}^{(r,k)}(x)+s_{n}^{(r,k)}(x)+nrB_{1}s_{n-1}^{(r,k)}(x)-r\sum_{l=0}^{n-2}(-1)^{n-l}\binom{n}{l}B_{n-l}s_{l}^{(r,k)}(x)\nonumber\\
&\quad -s_{n}^{(r+1,k)}(x)+s_{n}^{(r+1,k-1)}(x).\nonumber
\end{align}
Thus, by (\ref{eq:37}), we obtain the following theorem.

\begin{thm}\label{eq:thm3}
For $n\in\mathbf{N}$ with $n\geq 2$, we have
\begin{align*}
&ns_{n}^{(r,k)}(x)+n\left(\frac{1}{2}r-x\right)s_{n-1}^{(r,k)}(x)+r\sum_{l=0}^{n-2}(-1)^{n-l}\binom{n}{l}B_{n-l}s_{l}^{(r,k)}(x)\\
&=-s_{n}^{(r+1,k)}(x)+s_{n}^{(r+1,k-1)}(x).
\end{align*}
\end{thm}

\noindent For $r=0$, by Theorem \ref{eq:thm3}, we get
\begin{align*}
&nB_{n}^{(k)}(x)-nxB_{n-1}^{(k)}(x)\\
&=-B_{n}^{(k)}(x)+\frac{1}{2}nB_{n-1}^{(k)}(x)-\sum_{l=0}^{n-2}\binom{n}{l}B_{n-l}B_{l}^{(k)}(x)+\sum_{l=0}^{n}\binom{n}{l}B_{n-l}B_{l}^{(k-1)}(x).
\end{align*}
From (\ref{eq:4}), we note that
\begin{align}\label{eq:38}
s_{n}^{(r,k)}(y)&=\left\langle\left(\frac{t}{e^{t}-1}\right)^{r}\frac{Li_{k}\left(1-e^{-t}\right)}{1-e^{-t}}e^{yt}\Bigg\vert x^{n}\right\rangle\\
&=\left\langle\left(\frac{t}{e^{t}-1}\right)^{r}\frac{Li_{k}\left(1-e^{-t}\right)}{1-e^{-t}}e^{yt}\Bigg\vert xx^{n-1}\right\rangle\nonumber\\
&=\left\langle\partial_{t}\left(\left(\frac{t}{e^{t}-1}\right)^{r}\frac{Li_{k}\left(1-e^{-t}\right)}{1-e^{-t}}e^{yt}\right)\Bigg\vert x^{n-1}\right\rangle\nonumber\\
&=\left\langle\left(\partial_{t}\left(\frac{t}{e^{t}-1}\right)^{r}\right)\frac{Li_{k}\left(1-e^{-t}\right)}{1-e^{-t}}e^{yt}\Bigg\vert x^{n-1}\right\rangle\nonumber\\
&\quad+\left\langle\left(\frac{t}{e^{t}-1}\right)^{r}\left(\partial_{t}\frac{Li_{k}\left(1-e^{-t}\right)}{1-e^{-t}}\right)e^{yt}\Bigg\vert x^{n-1}\right\rangle\nonumber\\
&\quad+\left\langle\left(\frac{t}{e^{t}-1}\right)^{r}\frac{Li_{k}\left(1-e^{-t}\right)}{1-e^{-t}}\partial_{t}e^{yt}\Bigg\vert x^{n-1}\right\rangle.\nonumber
\end{align}
Therefore, by (\ref{eq:38}), we obtain the following theorem.

\begin{thm}\label{eq:thm4}
For $n\geq 1$, $r, k\in\mathbf{Z}$, we have
\begin{align*}
s_{n}^{(r,k)}(x)&=-rs_{n-1}^{(r,k)}(x)+r\sum_{l=0}^{n-1}\frac{\binom{n-1}{l}}{(n+1-l)(n-l)}s_{l}^{(r+1,k)}(x)\\
&\quad+\sum_{l=0}^{n-1}\left\{(-1)^{n-1-l}\binom{n-1}{l}\sum_{m=0}^{n-1-l}(-1)^{m}\frac{(m+1)!}{(m+2)^{k}}S_{2}(n-1-l,m)\right\}\\
&\quad\times\mathbb{B}_{l}^{(r)}(x-1)+xs_{n-1}^{(r,k)}(x).
\end{align*}
\end{thm}

\noindent Now, we compute
\begin{equation*}
\left\langle\left(\frac{t}{e^{t}-1}\right)^{r}Li_{k}\left(1-e^{-t}\right)\Bigg\vert x^{n+1}\right\rangle
\end{equation*}
in two different ways.\\
On the one hand,
\begin{align}\label{eq:39}
&\left\langle\left(\frac{t}{e^{t}-1}\right)^{r}Li_{k}\left(1-e^{-t}\right)\Bigg\vert x^{n+1}\right\rangle\\
&=\left\langle\left(\frac{t}{e^{t}-1}\right)^{r}\frac{Li_{k}\left(1-e^{-t}\right)}{1-e^{-t}}\Bigg\vert \left(1-e^{-t}\right)x^{n+1}\right\rangle\nonumber\\
&=\left\langle\left(\frac{t}{e^{t}-1}\right)^{r}\frac{Li_{k}\left(1-e^{-t}\right)}{1-e^{-t}}\Bigg\vert x^{n+1}-(x-1)^{n+1}\right\rangle\nonumber\\
&=\sum_{m=0}^{n}\binom{n+1}{m}(-1)^{n-m}\left\langle 1\big\vert s_{m}^{(r,k)}(x)\right\rangle\nonumber\\
&=\sum_{m=0}^{n}\binom{n+1}{m}(-1)^{n-m}s_{n}^{(r,k)}.\nonumber
\end{align}
On the other hand, we have
\begin{align}\label{eq:40}
&\left\langle\left(\frac{t}{e^{t}-1}\right)^{r}Li_{k}\left(1-e^{-t}\right)\Bigg\vert x^{n+1}\right\rangle=\left\langle Li_{k}\left(1-e^{-t}\right)\Bigg\vert \left(\frac{t}{e^{t}-1}\right)^{r}x^{n+1}\right\rangle\\
&=\left\langle Li_{k}\left(1-e^{-t}\right)\Bigg\vert \mathbb{B}_{n+1}^{(r)}(x)\right\rangle=\left\langle\int_{0}^{t}\left(Li_{k}\left(1-e^{-s}\right)\right)'ds\bigg\vert \mathbb{B}_{n+1}^{(r)}(x)\right\rangle\nonumber\\
&=\sum_{l=0}^{\infty}\sum_{m=0}^{l}\binom{l}{m}(-1)^{l-m}B_{n}^{(k-1)}\frac{1}{l!}\left\langle\int_{0}^{t}s^{l}ds\bigg\vert\mathbb{B}_{n+1}^{(r)}(x)\right\rangle\nonumber\\
&=\sum_{l=0}^{n}\sum_{m=0}^{l}\binom{l}{m}(-1)^{l-m}\frac{B_{m}^{(k-1)}}{(l+1)!}\left\langle 1\big\vert t^{l+1}\mathbb{B}_{n+1}^{(r)}(x)\right\rangle\nonumber\\
&=\sum_{l=0}^{n}\sum_{m=0}^{l}(-1)^{l-m}\binom{l}{m}\binom{n+1}{l+1}B_{m}^{(k-1)}\mathbb{B}_{n-l}^{(r)}.\nonumber
\end{align}
Therefore, by (\ref{eq:39}) and (\ref{eq:40}), we obtain the following theorem.

\begin{thm}\label{eq:thm5}
For $n\in\mathbf{Z}_{\geq 0}$, $r, k\in\mathbf{Z}$, we have
\begin{align*}
&\sum_{m=0}^{n}\binom{n+1}{m}(-1)^{n-m}s_{n}^{(r,k)}\\
&=\sum_{l=0}^{n}\sum_{m=0}^{l}(-1)^{l-m}\binom{l}{m}\binom{n+1}{l+1}B_{m}^{(k-1)}\mathbb{B}_{n-l}^{(r)}.
\end{align*}
\end{thm}

\begin{lem}[$\lbrack 5\rbrack$]\label{eq:lem6}
For $k\in\mathbf{Z}$ and $m\geq 1$, we have
\begin{align}\label{eq:41}
&\left(\sum_{l=0}^{m}\left[ \begin{array}{c}
m\\
l
\end{array} \right]\partial_{t}^{l}\right)\frac{Li_{k}\left(1-e^{-t}\right)}{1-e^{-t}}\\
&=\frac{1}{\left(e^{t}-1\right)^{m}}\sum_{l=0}^{m}(-1)^{m-l}\left[ \begin{array}{c}
m+1\\
l+1
\end{array} \right]\frac{Li_{k-l}\left(1-e^{-t}\right)}{1-e^{-t}},\nonumber
\end{align}
where $\left[ \begin{array}{c}
m\\
l
\end{array} \right]=\lvert S_{1}(m,l)\rvert$ and $S_{1}(m,l)$ is the stirling number of the first kind.
\end{lem}

\noindent Now, we compute $\left\langle\left(\frac{t}{e^{t}-1}\right)^{m}\sum_{l=0}^{m}(-1)^{m-l}\left[ \begin{array}{c}
m+1\\
l+1
\end{array} \right]\frac{Li_{k-l}\left(1-e^{-t}\right)}{1-e^{-t}}\bigg\vert x^{n}\right\rangle$ in two different ways.\\
On the one hand,
\begin{align}\label{eq:42}
&\left\langle \left(\frac{t}{e^{t}-1}\right)^{m}\sum_{l=0}^{m}(-1)^{m-l}\left[ \begin{array}{c}
m+1\\
l+1
\end{array} \right]\frac{Li_{k-l}\left(1-e^{-t}\right)}{1-e^{-t}}\Bigg\vert x^{n}\right\rangle\\
&=\left\langle 1\Bigg\vert \sum_{l=0}^{m}(-1)^{m-l}\left[ \begin{array}{c}
m+1\\
l+1
\end{array} \right]\left(\frac{t}{e^{t}-1}\right)^{m}\frac{Li_{k-l}\left(1-e^{-t}\right)}{1-e^{-t}}x^{n}\right\rangle\nonumber\\
&=\left\langle 1\Bigg\vert\sum_{l=0}^{m}(-1)^{m-l}\left[ \begin{array}{c}
m+1\\
l+1
\end{array} \right]s_{n}^{(m,k-l)}(x)\right\rangle\nonumber\\
&=\sum_{l=0}^{m}(-1)^{m-l}\left[ \begin{array}{c}
m+1\\
l+1
\end{array} \right]s_{n}^{(m,k-l)}\nonumber
\end{align}
On the other hand, by Lemma \ref{eq:lem6}, it is equal to 
\begin{align}\label{eq:43}
&\left\langle\left(\frac{t}{e^{t}-1}\right)^{m}\sum_{l=0}^{m}(-1)^{m-l}\left[ \begin{array}{c}
m+1\\
l+1
\end{array} \right]\frac{Li_{k-l}\left(1-e^{-t}\right)}{1-e^{-t}}\Bigg\vert x^{n}\right\rangle\\
&=\left\langle t^{m}\left(\sum_{l=0}^{m}\left[ \begin{array}{c}
m\\
l
\end{array} \right]\partial_{t}^{l}\right)\frac{Li_{k}\left(1-e^{-t}\right)}{1-e^{-t}}\Bigg\vert x^{n}\right\rangle\nonumber\\
&=\left\langle\sum_{l=0}^{m}\left[ \begin{array}{c}
m\\
l
\end{array} \right]\partial_{t}^{l}\frac{Li_{k}\left(1-e^{-t}\right)}{1-e^{-t}}\Bigg\vert t^{m}x^{n}\right\rangle\nonumber\\
&=\left\{ \begin{array}{cc}
(n)_{m}\left\langle\sum_{l=0}^{m}\left[ \begin{array}{c}
m\\
l
\end{array} \right]\partial_{t}^{l}\frac{Li_{k}\left(1-e^{-t}\right)}{1-e^{-t}}\Bigg\vert x^{n-m}\right\rangle, & \text{if} \,\,n\geq m,\\
0, & \text{if}\,\, 0\leq n\leq m-1.
\end{array} \right.\nonumber
\end{align}
For $n\geq m$, we have
\begin{align}\label{eq:44}
&(n)_{m}\left\langle\sum_{l=0}^{m}\left[\begin{array}{c}
m\\
l
\end{array}\right]
\partial_{t}^{l}\frac{Li_{k}\left(1-e^{-t}\right)}{1-e^{-t}}\Bigg\vert x^{n-m}\right\rangle\\
&=(n)_{m}\left\langle\frac{Li_{k}\left(1-e^{-t}\right)}{1-e^{-t}}\Bigg\vert\left(\sum_{l=0}^{m}\left[ \begin{array}{c}
m\\
l
\end{array} \right]x^{l}\right)x^{n-m}\right\rangle\nonumber\\
&=(n)_{m}\sum_{l=0}^{m}\left[ \begin{array}{c}
m\\
l
\end{array} \right]\left\langle 1\Bigg\vert \frac{Li_{k}\left(1-e^{-t}\right)}{1-e^{-t}}x^{n-m+l}\right\rangle\nonumber\\
&=(n)_{m}\sum_{l=0}^{m}\left[ \begin{array}{c}
m\\
l
\end{array} \right]B_{n-m+l}^{(k)}.\nonumber
\end{align}
Therefore, by (\ref{eq:42}), (\ref{eq:43}) and (\ref{eq:44}), we obtain the following theorem.

\begin{thm}[$\lbrack 5\rbrack$]\label{eq:thm7}
For $k\in\mathbf{Z}$, $m\geq 1$, we have
\begin{align*}
&\sum_{l=0}^{m}(-1)^{m-l}\left[ \begin{array}{c}
m+1\\
l+1
\end{array} \right]s_{n}^{(m,k-l)}\\
&=\left\{ \begin{array}{cc}
(n)_{m}\sum_{l=0}^{m}\left[ \begin{array}{c}
m\\
l
\end{array} \right]B _{n-m+l}^{(k)}, & \text{if}\,\, n\geq m,\\
0, & \text{if}\,\, 0\leq n\leq m-1.
\end{array} \right.
\end{align*}
\end{thm}

\noindent Now, we consider the following two Sheffer sequences:
\begin{align}\label{eq:45}
s_{n}^{(r,k)}(x)\sim\left(\left(\frac{e^{t}-1}{t}\right)^{r}\frac{1-e^{-t}}{Li_{k}\left(1-e^{-t}\right)}, t\right)
\end{align}
and
\begin{align*}
E_{n}
^{(s)}(x)\sim\left(\left(\frac{e^{t}+1}{2}\right)^{s},t\right).
\end{align*}
Let us assume that
\begin{equation}\label{eq:46}
s_{n}^{(r,k)}(x)=\sum_{m=0}^{n}C_{n,m}E_{m}^{(s)}(x).
\end{equation}
Then, from (\ref{eq:18}), we have
\begin{align}\label{eq:47}
C_{n,m}&=\frac{1}{m!}\left\langle\frac{\left(\frac{e^{t}+1}{2}\right)^{s}}{\left(\frac{e^{t}-1}{t}\right)^{r}\frac{1-e^{-t}}{Li_{k}\left(1-e^{-t}\right)}}t^{m}\Bigg\vert x^{n}\right\rangle\\
&=\frac{1}{m!}\left\langle\left(\frac{e^{t}+1}{2}\right)^{s}\left(\frac{t}{e^{t}-1}\right)^{r}\frac{Li_{k}\left(1-e^{-t}\right)}{1-e^{-t}}\Bigg\vert t^{m}x^{n}\right\rangle\nonumber\\
&=\frac{\binom{n}{m}}{2^{s}}\sum_{j=0}^{s}\binom{s}{j}\left\langle e^{jt}\Bigg\vert\left(\frac{t}{e^{t}-1}\right)^{r}\frac{Li_{k}\left(1-e^{-t}\right)}{1-e^{-t}}x^{n-m}\right\rangle\nonumber\\
&=\frac{\binom{n}{m}}{2^{s}}\sum_{j=0}^{s}\binom{s}{j}s_{n-m}^{(r,k)}(j).\nonumber
\end{align}
Therefore, by (\ref{eq:46}) and (\ref{eq:47}), we obtain the following theorem.

\begin{thm}\label{eq:thm8}
For $r, k\in\mathbf{Z}$, $n, m\in\mathbf{Z}_{\geq 0}$, we have
\begin{align*}
s_{n}^{(r,k)}(x)=\frac{1}{2^{s}}\sum_{m=0}^{n}\left\{\binom{n}{m}\sum_{j=0}^{s}\binom{s}{j}s_{n-m}^{(r,k)}(j)\right\}E_{m}^{(s)}(x).
\end{align*}
\end{thm}

\noindent Let us consider the following two Sheffer sequences:
\begin{equation}\label{eq:48}
s_{n}^{(r,k)}(x)\sim\left(\left(\frac{e^{t}-1}{t}\right)^{r}\frac{1-e^{-t}}{Li_{k}\left(1-e^{-t}\right)}, t\right),\quad \mathbb{B}_{n}^{(s)}(x)\sim\left(\left(\frac{e^{t}-1}{t}\right)^{s}, t\right).
\end{equation}
Let
\begin{equation}\label{eq:49}
s_{n}^{(r,k)}(x)=\sum_{m=0}^{n}C_{n,m}\mathbb{B}_{m}^{(s)}(x).
\end{equation}
From (\ref{eq:18}), we note that
\begin{align}\label{eq:50}
C_{n,m}&=\frac{1}{m!}\left\langle\left(\frac{t}{e^{t}-1}\right)^{r-s}\frac{Li_{k}\left(1-e^{-t}\right)}{1-e^{-t}}t^{m}\Bigg\vert x^{n}\right\rangle\\
&=\binom{n}{m}\left\langle 1\Bigg\vert\left(\frac{t}{e^{t}-1}\right)^{r-s}\frac{Li_{k}\left(1-e^{-t}\right)}{1-e^{-t}}x^{n-m}\right\rangle\nonumber\\
&=\binom{n}{m}s_{n-m}^{(r-s,k)}.\nonumber
\end{align}
Therefore, by (\ref{eq:49}) and (\ref{eq:50}), we obtain the following theorem.

\begin{thm}\label{eq:thm9}
For $r, s\in\mathbf{Z}$, $n\in\mathbf{Z}_{\geq 0}$, we have
\begin{equation*}
s_{n}^{(r,k)}(x)=\sum_{m=0}^{n}\binom{n}{m}s_{n-m}^{(r-s,k)}\mathbb{B}_{m}^{(s)}(x).
\end{equation*}
\end{thm}

\noindent We note that
\begin{equation*}
s_{n}^{(r,k)}(x)\sim\left(\left(\frac{e^{t}-1}{t}\right)^{r}\frac{1-e^{-t}}{Li_{k}\left(1-e^{-t}\right)}, t\right),
\end{equation*}
and
\begin{equation*}
H_{n}^{(s)}(x\vert\lambda)\sim\left(\left(\frac{e^{t}-\lambda}{1-\lambda}\right)^{s}, t\right),\,\,\,\,(\text{see}\,\,\lbrack 1,6\rbrack).
\end{equation*}
By the same method, we get
\begin{equation}\label{eq:51}
s_{n}^{(r,k)}(x)=\frac{1}{\left(1-\lambda\right)^{s}}\sum_{m=0}^{n}\left\{\binom{n}{m}\sum_{j=0}^{s}\binom{s}{j}\left(-\lambda\right)^{s-j}s_{n-m}^{(r,k)}(j)\right\}H_{m}^{(s)}(x\vert\lambda).
\end{equation}


\bigskip
ACKNOWLEDGEMENTS. This work was supported by the National Research Foundation of Korea(NRF) grant funded by the Korea government(MOE) (No.2012R1A1A2003786 ).
\bigskip


\noindent
\author{Department of Mathematics, Sogang University, Seoul 121-742, Republic of Korea
\\e-mail: dskim@sogang.ac.kr}\\
\\
\author{Department of Mathematics, Kwangwoon University, Seoul 139-701, Republic of Korea
\\e-mail: tkkim@kw.ac.kr}
\end{document}